\documentclass[11pt]{amsart}
\usepackage{amscd,amssymb}
\usepackage{hyperref} 

\input xy
\xyoption{all}

\newcommand{\hsk}{\hskip 0.3cm}
\newcommand{\ssk}{\vskip 0.2cm}
\newcommand{\sk}{\vskip 0.4cm}
\newcommand{\bsk}{\vskip 1cm}

\newtheorem{thm}{Theorem}[section]
\newtheorem{corol}[thm]{Corollary}
\newtheorem{lemma}[thm]{Lemma}
\newtheorem{prop}[thm]{Proposition}

\theoremstyle{definition}
\newtheorem{defin}[thm]{Definition}

\theoremstyle{remark}
\newtheorem{remark}[thm]{Remark}

\newtheorem{example}[thm]{Example}
\newtheorem{examples}[thm]{Examples}
\numberwithin{equation}{section} 

\def\R {{\Bbb R}}

\def\N{{\Bbb N}}

\newcommand{\Acal}{\mathcal{A}}
\newcommand{\Fcal}{\mathcal{F}}

\newcommand{\ben}{\begin{enumerate}}

\newcommand{\een}{\end{enumerate}}
\newcommand{\bit}{\begin{itemize}}

\newcommand{\eit}{\end{itemize}}

\def\a{\alpha}

\def\eps{\varepsilon}

\def\QED{\nobreak\quad\ifmmode\roman{Q.E.D.}\else{\rm Q.E.D.}\fi}

\begin{document}

%

\title[]
{Compactifications of Semigroups and Semigroup Actions}

\author[]{Michael Megrelishvili}
\address{Department of Mathematics,
Bar-Ilan University, 52900 Ramat-Gan, Israel}
\email{megereli@math.biu.ac.il}
\urladdr{http://www.math.biu.ac.il/$^\sim$megereli}

\date{December 7, 2006}


\begin{abstract}
An action of a topological semigroup $S$ on $X$ is {\it
compactifiable} if this action is a restriction of a jointly
continuous action of $S$ on a Hausdorff compact space $Y$. A
topological semigroup S is compactifiable if the left action of
$S$ on itself is compactifiable. It is well known that every
Hausdorff topological group is compactifiable. This result cannot
be extended to the class of Tychonoff topological monoids. At the
same time, several natural constructions lead to compactifiable
semigroups and actions.


  We prove that the semigroup $C(K,K)$ of all continuous selfmaps on the
Hilbert cube $K=[0,1]^{\omega}$ is a universal second countable
compactifiable semigroup ({\it semigroup version of Uspenskij's
theorem}). Moreover, the Hilbert cube $K$ under the action of
$C(K,K)$ is universal in the realm of all compactifiable $S$-flows
$X$ with compactifiable $S$ where both $X$ and $S$ are second
countable.

We strengthen some related results of Kocak \& Strauss \cite{KS}
and Ferry \& Strauss \cite{FS} about Samuel compactifications of
semigroups. Some results concern compactifications with separately
continuous actions, LMC-compactifications and LMC-functions
introduced by Mitchell.
\end{abstract}


\maketitle





 \tableofcontents

\section{Introduction}
\sk

A major role of semigroup actions and semigroup compactifications
is now well understood. See for example the books \cite{BJMo,BJM} and
\cite{Ru}.
Very little is known however about sufficient conditions which
ensure the existence of {\it proper} compactifications in the case
of monoidal actions. This contrasts the case of topological group
actions (see for example \cite{vr-can, vrexist, vr-embed,
lud-vr,MeSc1, me-singap, Me1, Me2, Me-opit2}).

A semigroup action $S \times X \to X$, or, a {\it flow} $(S, X)$,
is {\it compactifiable} if there exists a {\it proper}
$S$-compactification $X \hookrightarrow Y$. That is, if the
original action is a restriction of a jointly continuous action on
a Hausdorff compact $S$-flow $Y$. In this article we require that
$S$ is a {\it topological semigroup} (the multiplication is
jointly continuous).  We say that a topological semigroup $S$ is
{\it compactifiable} if the flow $(S,S)$, the regular left action,
is compactifiable. Passing to the {\it Ellis semigroup} $E(Y)$ of
an $S$-compactification $Y$ of a monoid $S$ we see that $S$ is
compactifiable iff $S$ has a proper {\it dynamical
compactification} in the sense of Ruppert \cite{Ru} (see also the
{\it monoidal compactification} in the sense of Lawson \cite{La}).

If a topological semigroup $S$ algebraically is a group we say
that $S$ is a {\it paratopological group}. As usual, {\it
topological group} means that in addition we require the
continuity of the inverse operation. Due to Teleman \cite{Te}
every Hausdorff (equivalently: Tychonoff) topological group is
compactifiable. This classical result cannot be extended to the
class of Tychonoff topological semigroups.
For instance, the multiplicative monoid $([0,\infty), \cdot)$ of
all nonnegative reals is not compactifiable (see Example
\ref{ex:count}.2 below) and even not {\it LMC-compactifiable} as
it follows by a result of Hindman and Milnes \cite{HM}. The latter
means in fact
%
that there is no proper $S$-compactification $S \to Y$ with a
separately continuous action on $Y$. LMC is an abbreviation of
\emph{Left Multiplicatively Continuous}.  LMC-compactifications
and LMC-functions for semigroups were introduced by Mitchell,
\cite{Mi,HM,BJMo}. The case of separately continuous
compactifications is parallel to the theory of right topological
compactifications
 and \emph{generalized $LMC$-functions} (see Definition \ref{d:LMC}).
This direction is linked to Banach representations of semigroups
and actions (in the sense of \cite{me-nz}) and to corresponding
generalized matrix coefficients.

 One of our aims in the
present paper is to study the similarities and differences in the theory of flow compactifications
when we pass from groups to semigroups. We emphasize the
limitations providing several non-compactifiable semigroups and actions
with ``good topological properties'' (contrasting the case of topological groups).

The classical Gelfand-Naimark 1-1 correspondence between Banach
subalgebras of $C(X)$ and the compactifications of $X$ can be
extended to the category of $S$-flows describing jointly
continuous $S$-compactifications by subalgebras of the algebra
$RUC_S(X)$ of all {\it right uniformly continuous functions} on
$X$ (see Definition \ref{d:ruc}).
This theory is well known for topological group actions
(see, for example, J. de Vries \cite{vr-embed}).
One can easily extend it to the case of topological semigroup actions.
Some results in this direction can be found in the work of Ball and Hagler \cite{BH}.

We establish some sufficient and necessary conditions in terms of
uniform structures. In particular, we strengthen two results of
Kocak and Strauss \cite{KS} and also a result of Ferry and Strauss
\cite{FS} (see Corollary \ref{c:str} and Remark
\ref{r:results}.1).

The topological monoid $C(K,K)$
of all continuous self-maps endowed with the compact open topology
is compactifiable.
If $E$ is a normed space then the monoid $(\Theta(E), norm)$
of all contractive linear self-operators $E \to E$ is
compactifiable endowed with the norm topology. It is not true with
respect to the strong operator topology $\tau_{s}$ on $\Theta(E)$.
However, its topological opposite semigroup
$(\Theta(E)^{op},\tau_{s})$ is always compactifiable.

%

A paratopological group $G$ is compactifiable iff $G$ is a
topological group. It follows in particular, that the Sorgenfrey
Line, as an additive monoid, is not compactifiable.

One of our main results states that the semigroup
$U:=C(I^{\omega}, I^{\omega})$ is a universal second countable
compactifiable semigroup. It is a semigroup version of Uspenskij's
theorem \cite{Us1} about universality of the group
$Homeo(I^{\omega})$. Moreover, strengthening a result of
\cite{Me2}, we establish that the action of $U$ on $I^{\omega}$ is
universal in the realm of compactifiable $S$-flows $X$ (with
compactifiable $S$) where $X$ and $S$ both are separable and
metrizable.

The present paper influenced especially by \cite{FS, KS, Pe1,
Us1}.

\vskip 1cm
\section{Semigroup actions: natural examples and representations}
\vskip 0.4cm


Let $\pi: P \times X \to Z$ be a map. For $p_0 \in P$ and $x_0 \in
X$ define \emph{left and right translations} by
$$\lambda_{p_0}: X \to Z, \hsk x \mapsto \pi(p_0,x)$$ and
$$\rho_{x_0}: P \to Z, \hsk p \mapsto \pi(p,x_0)$$ respectively.
The map $\pi$ is \emph{left (right) continuous} if every left
(right) translation is continuous.

\begin{lemma} \label{l:contpoints}
Let $\pi: P \times X \to Z$ be a right continuous map, $P'$ and
$X'$ be dense subsets of $P$ and $X$ respectively. Assume that the
map $\lambda_{p'} : X \to Z$ is continuous for every $p' \in P'$
and $Z$ is a regular space. Then if $P' \times X' \to Z$ is
continuous at $(p',x')$ then $\pi: P \times X \to Z$ is continuous
at $(p',x')$.
\end{lemma}
\begin{proof}
Let $O$ be a neighborhood of $\pi(p'x')$ in $Z$. Since $Z$ is
regular one can choose a neighborhood $U$ of $\pi(p'x')$ such that
$cl (U) \subset O$. Now by continuity of $\pi'$ at $(p',x')$
choose the neighborhoods $V$ of $p'$ in $P'$ and $W$ of $x'$ in
$X'$ s.t. $\pi'(t,y) \in U$ for every $(t,y) \in V \times W$. Now
$\pi(p,x) \in cl(U) \subset O$ for every $p \in cl(V)$ and $x \in
cl(W)$. Indeed choose two nets $a_i \in V$ and $b_j \in W$ s.t.
$\lim_ia_i=p$ in $P$ and $\lim_j b_j=x$ in $X$.

Since $a_i \in P'$ the map $\lambda_{a_i}$ is continuous for every
$i$. We have $\lim_j \pi(a_i,b_j) = \pi(a_i,x) \in cl(U)$ for
every $i$. Now by right continuity of $\pi$ we obtain $\lim_i
\pi(a_i,x)=\pi(p,x) \in cl(U)$. This implies the continuity of
$\pi$ at $(p',x')$ because $cl(V)$ and $cl(W)$ are neighborhoods
of $p'$ and $x'$ in $P$ and $X$ respectively.
\end{proof}

A topologized semigroup $S$ is: (a) \emph{left} (\emph{right})
topological; (b) \emph{semitopological}; (c) \emph{topological} if
the multiplication function $S\times S \to S$ is left (right)
continuous, separately continuous, or jointly continuous,
respectively.

A {\it topological (left) $S$-flow} (or an {\it $S$-space}) is a triple $(S,X,\pi)$
where $\pi : S\times X \to X$ is a jointly continuous left action of a topological
semigroup $S$ on a topological space $X$; we write it also as a pair $(S,X)$, or simply, $X$
(when $\pi$ and $S$ are understood).
As usual we write $sx$ instead of $\pi(s,x)={\breve s} (x)=
{\tilde x}(s)$. ``Action'' means that always
$s_1(s_2x)=(s_1s_2)x$. Every $x \in X$ defines the orbit map
${\tilde x}: S \to X, \ s \mapsto sx$. Every $s \in S$ gives rise
to the $s$-translation ${\breve s}: X \to X, \ x \mapsto sx$. The
action is {\it monoidal} If $S$ is a monoid and the identity $e$
of $S$ acts as the identity transformation of $X$.

If the action $S \times X \to X$ is separately continuous (that
is, all orbit maps ${\tilde x}$ and all translations ${\breve s}:
X \to X$ are continuous) then we say that $X$ (or, $(S,X)$) is a
{\it semitopological $S$-flow}.

 A {\it right flow} $(X,S)$ can be defined analogously. If
$S^{op}$ is the {\it opposite semigroup} of $S$ with the same
topology then $(X,S)$ can be treated as a left flow $(S^{op},X)$
(and vice versa).

Let $h: S_1 \to S_2$ be a semigroup homomorphism, $S_1$ act on
$X_1$ and $S_2$ on $X_2$. A map $\a: X_1 \to X_2$ is said to be
$h$-{\it equivariant} if $\a(sx)=h(s)\a(x)$ for every $(s,x) \in
S_1 \times X_1$. Sometimes we say that the pair $(h, \a)$ is {\it
equivariant}. For $S_1=S_2$ with $h=id_S$, we say: {\it $S$-map}.
The map $h: S_1 \to S_2$ is a {\it co-homomorphism} iff $S_1 \to
S_2^{op}, s \mapsto h(s)$ is a homomorphism. We say that $(h,\a)$
is {\it proper} if $\a$ is a topological embedding.

Let $\mu$ be a uniform structure on a set $X$. We assume that it
is separated. Then the induced topology $top(\mu)$ on $X$ is
Tychonoff. A uniformity $\mu$ on a topological space $(X,\tau)$ is
said to be {\it compatible} if $top(\mu)=\tau$. ``Compact'' will
mean compact and Hausdorff.

\sk
Recall some natural ways getting topological monoids and monoidal actions.

Let $V$ be a normed space. The closed unit ball of $V$ we denote
by $B_V$. Weak star compact unit ball $B_{V^*}$ in the dual space
$V^*$ will be denoted also by $B^*$.


\begin{examples}
\label{ex:monoids}
\ben
\item Let $(Y,\mu)$ be a uniform space.
Denote by $\mu_{sup}$ the uniformity  of uniform convergence on
the set $Unif(Y,Y)$ of all uniform self-maps of $Y$. Then under
the corresponding topology $top(\mu_{sup})$ on $Unif(Y,Y)$ and the
usual composition we get a topological monoid. For every
subsemigroup $S \subset Unif(Y,Y)$ the induced action $S \times Y
\to Y$ defines a topological flow.
\item For instance, for every compact space $Y$ the semigroup $C(Y,Y)$
 endowed with the compact open topology is a topological monoid. Note also
that the subset $Homeo(Y)$ in $C(Y,Y)$ of all homeomorphisms $Y
\to Y$ is a topological group.
\item
For every metric space $(M,d)$ the semigroup $\Theta(M,d)$ of all
{\it $d$-contractive maps} $f: X\to X$ (that is, $d(f(x),f(y))
\leq d(x,y)$) is a topological monoid with respect to the topology
of pointwise convergence. Furthermore, the map $\Theta(M,d) \times
M \to M$ is a jointly continuous monoidal action.
\item
For every normed space $(V, ||\cdot||)$ the semigroup $\Theta(V)$
of all contractive linear operators $V \to V$ endowed with the
\emph{strong operator topology} (being a topological submonoid of
$\Theta(V,d)$ where $d(x,y):=||x-y||$) is a topological monoid.
The subspace $Is(V)$ of all linear onto isometries is a
topological group.
\item For every normed space $V$ and a subsemigroup $S \subset \Theta(V)^{op}$
the induced action $S \times
B^* \to B^*$ on the compact space $B^*$ is jointly continuous (see Lemma \ref{l:opp}).
\item
Every normed algebra $A$ treated as a multiplicative monoid is a
topological monoid. The subset $B_A$ is a topological submonoid.
In particular, for every normed space $V$ the monoids $L(V)$ and
$B_{L(V)}$ of all bounded and, respectively, of all contractive
linear operators $V \to V$ are topological monoids endowed with
the norm topology. Observe that $B_{L(V)}$ and $\Theta(V)$ {\it
algebraically} are the same monoids.
 \een

\ssk
We omit the straightforward arguments.
\end{examples}

An action $S \times X \to X$ on a metric space $(X,d)$ is
\emph{contractive}
if every $s$-translation $\tilde{s}: X \to X$ lies in
$\Theta(X,d)$. It defines a natural homomorphism $h: S \to
\Theta(X,d)$.

\begin{remark} \label{r:simple}
\ben
\item If an action of $S$ on $(X.d)$ is contractive
then it is easy to show that the following
conditions are equivalent: \bit
\item [(i)] The action is jointly continuous.

\item [(ii)] The action is separately continuous.

\item [(iii)] The restriction $S \times Y \to X$ to some
dense subspace Y of X is separately continuous.

\item [(iv)] The natural homomorphism $h: S \to \Theta(X,d)$ is continuous.
\eit

\item If $j: V \hookrightarrow \widehat{V}$ is the completion of a
normed space $V$ then we have the following canonical equivariant
inclusion of monoidal actions
$$
(\Theta(V),V) \rightrightarrows (\Theta(\widehat{V}),
\widehat{V}).
$$
\een
\end{remark}

\sk

The Banach algebra of all continuous real valued  bounded
functions on a topological space $X$ will be denoted by $C(X)$.
Every left action $\pi: S \times X \to X$ induces the
co-homomorphism $h_{\pi}: S \to C(X)$ and the right action
$C(X)\times S \to C(X)$ where $(fs)(x)=f(sx)$. While the
translations are continuous, the orbit maps ${\tilde f}: S \to
C(X)$ are not necessarily norm (even weakly) continuous and
requires additional assumptions (see Definition \ref{d:ruc}).


For every normed space $V$ the usual adjoint map $adj: L(V) \to
L(V^*), \hsk \phi \mapsto \phi^*$ is an injective co-homomorphism
of monoids.

The following two simple lemmas are very useful. For some closely
related results see \cite{Us1}, \cite[Chapter 5]{Ak} and
\cite[Fact 2.2]{me-nz} .

\begin{lemma} \label{l:opp}
For every normed space $V$ the injective map $\gamma:
\Theta(V)^{op} \hookrightarrow C(B^*,B^*)$, induced by the adjoint
map $adj: L(V) \to L(V^*)$, is a topological (even uniform) monoid
embedding. In particular,
$$
\Theta(V)^{op} \times B^* \to B^*
$$
is a jointly continuous monoidal action of $\Theta(V)^{op}$ on the
compact space $B^*$.
\end{lemma}
\begin{proof}
The strong uniformity on $\Theta (V)$ is generated by the family
of pseudometrics $\{p_v: v \in V\}$,where $p_v(s,t)=||sv-tv||$. On
the other hand the family of pseudometrics $\{q_v: v \in
V\}$,where $q_v(s,t)=sup\{(fs)(v)-(ft)(v): f\in B^*\}$ generates
the natural uniformity inherited from $C(B^*,B^*)$. Now observe that
$p_v(s,t)=q_v(s,t)$. This proves that $\gamma$ is a uniform (and
hence, also, topological) embedding.
%
\end{proof}


\begin{lemma} \label{Teleman}
Let $V$ be a Banach space. Suppose that $\pi: V \times S \to V$ is
a right action of a topologized semigroup $S$ by linear
contractive operators. The following are equivalent: \bit
\item [(i)]
The co-homomorphism $h: S \to \Theta(V), \ h(s)(v):=vs$ is
strongly continuous.
\item [(ii)]
 The induced
action $S \times B^* \to B^*, \ (s\psi)(v):=\psi(vs)$ on the weak
star compact ball $B^*$ is jointly continuous. \eit
\end{lemma}
\begin{proof} (i) $\Rightarrow$ (ii): Let $h: S \to \Theta(V)$ be strongly continuous. Then
by Lemma \ref{l:opp} the composition $\gamma \circ h: S \to
C(B^*,B^*)$ is also continuous. This yields (ii) (see Example
\ref{ex:monoids}.2).

(i) $\Leftarrow$ (ii): Since the action $S \times B^* \to B^*$ is
continuous and $B^*$ is compact the homomorphism $S \to
C(B^*,B^*), \ s \mapsto \breve{s}$ is continuous. Again by Lemma
\ref{l:opp} we get that the co-homomorphism $h: S \to \Theta(V)$
is strongly continuous.
%
\end{proof}

\begin{defin} \label{d:repr}
\ben
\item
\cite[Definition 3.1]{me-nz} \ A (continuous) {\it representation
of a flow $(S,X)$ on a normed space} $V$ is an equivariant pair
$$
(h, \alpha): (S,X) \rightrightarrows (\Theta(V)^{op}, B^*)
$$
where $\alpha: X \to B^*$ is weak$^*$ continuous and $h: S \to
\Theta(V)^{op}$ is a (resp.: strongly continuous) homomorphism.


\item
A {\it representation of $(S,X)$ on a uniform space} $(Y,\mu)$ is an
equivariant pair
$$
(h, \alpha): (S,X) \rightrightarrows (Unif(Y,Y), Y)
$$
where $h: S \to Unif(Y,Y)$ is a continuous homomorphism and $\a: X
\to (Y,top(\mu))$ is a continuous map (cf. Definition
\ref{d:equiun}.3).

\een
\end{defin}

\begin{defin} \label{d:equic}
\ben
\item Let $S \times X \to X$ be a semigroup action.
A uniformity $\mu$ on $X$ is {\it equicontinuous} if for every
$\varepsilon \in \mu$ and any $x_0 \in X$ there exists a
neighborhood $O$ of $x_0$ such that $(sx,sx_0) \in \varepsilon$
for every $x \in O$ and every $s \in S$. If there exists $\delta
\in \mu$ such that $(sx,sy) \in \varepsilon$ holds for every pair
$x,y$ from $X$ then as usual we say that $\mu$ is {\it uniformly
equicontinuous}. In the case of right actions the definitions are
similar.
\item A pseudometric $d$
on a semigroup $S$ is {\it right contractive} if $d(xs,ys) \leq
d(x,y)$ for every $x,y,s \in S$.
\item
A uniform structure $\mu$ on a semigroup $S$ is \emph{right
invariant} (see also \cite[p. 98]{FS} and Lemma \ref{l:rightinv})
if for every $\eps \in \mu$ there exists $\delta \in \mu$ such
that $\delta \subset \eps$ and $(sx,tx) \in \delta$ for every
$(s,t) \in \delta$, $x \in S$.
\een
\end{defin}

\begin{lemma} \label{l:rightinv} Let $\mu$ be a compatible uniform structure
on a topological semigroup $S$. The following conditions are
equivalent: \ben
\item
 $\mu$ can be generated by a
family of right contractive pseudometrics.
\item
$\mu$ is right invariant on $S$.
\item
The right action of $S$ on itself is $\mu$-uniformly
equicontinuous (that is, for every $\varepsilon \in \mu$ there
exists $\delta \in \mu$ such that $(sx,tx) \in \varepsilon$ for
every $(s,t) \in \delta$, $x \in S$). \een
\end{lemma}
\begin{proof}
The implications (1) $\Rightarrow$ (2) and (2) $\Rightarrow$ (3)
are trivial.

(3) $\Rightarrow$ (1): Assume that the right action of $S$ on
itself is $\mu$-uniformly equicontinuous. Choose a family
$\{d_i\}_{i \in I}$ of pseudometrics on S which generates the
uniformity $\mu$. For every $i \in I$ define $$d^*_i(x,y):=max
\{sup_{s \in S} d_i(xs,ys), d(x,y) \}$$ Then the new system
$\{d^*_i\}_{i \in I}$ consists by right contractive pseudometrics
and still generates the same uniformity $\mu$.
\end{proof}

\begin{example} \label{e:unif-eq}
\ben
\item For every topological group $G$ the
right uniformity $\mathcal R(G)$ of $G$ is the \emph{unique} right
invariant compatible uniformity on $G$, \cite[Lemma 2.2.1]{RD}.
\item
 Let $(X,\mu)$ be a uniform space and $\mu_{sup}$ be the
corresponding natural uniformity on $Unif(X,X)$. Assume that $S$
is a subsemigroup of $Unif(X,X)$. Then the subspace uniformity
$\mu_{sup}|_S$ on $S$ is right invariant.
 \een
\end{example}

The following proposition is an equivariant version of the well
known Arens-Eells embedding construction \cite{AE}.

\begin{prop} \label{l:AE}
Let $S \times X \to X$ be a continuous contractive action of a
semigroup $S$ on a
bounded metric space $(X,d)$.
Then there exists a normed (equivalently: {\it Banach}) space $E$
and an equivariant pair
$$
(h,\a): (S,X) \rightrightarrows (\Theta(E),E)
$$
such that $h: S \to \Theta(E)$ is a strongly continuous homomorphism and
$\a: X \to E$ is an isometric embedding.
\end{prop}
\begin{proof}
By Remark \ref{r:simple} it suffices to
give a proof for normed $E$.
Since the metric is bounded we can suppose that $X$ contains a
fixed point $z$ (adjoining if necessary a fixed point $z$ and
defining $d(x,z)=diam (X,d) <\infty$ for every $x \in X$). We can
use the {\it Arens-Eells isometric embedding}
$$i: X \to A(X), \ x \mapsto x-z$$
(see \cite{AE}) of a pointed metric space $(X,z,d)$ into a normed
space $(A(X), ||\cdot ||)$.
The elements of $A(X)$ are the formal sums of the form
$\sum_{i=1}^n c_i(x_i-y_i)$, where $x_i, y_i \in X$ and $c_i \in
\R$. Define the natural left action
$$S \times A(X) \to A(X), \ \
 \ s\sum_{i=1}^n c_i(x_i-y_i):=\sum_{i=1}^n
c_i(sx_i-sy_i).$$
The desired norm on $A(X)$ is defined by setting
$$||u||:=inf \sum_{i=1}^n |c_i|d(x_i,y_i)),$$
where we compute the infimum with respect to the all presentations
of $u \in A(X)$ as the sums $u=\sum_{i=1}^n c_i(x_i-y_i)$ with
$x_i, y_i \in X$. This explicit description shows that $||su||
\leq ||u||$ for every $s \in S$ because $d(sx_i,sy_i) \leq
d(x_i,y_i)$. Therefore the action $S \times X \to X$ can be
extended to the canonically defined action $S \times A(X) \to
A(X)$ by contractive linear operators. Moreover it is clear that
every orbit mapping $S \to A(X), \ s \mapsto su$ is continuous for
every $u \in A(X)$. Thus we get a continuous homomorphism $h: S
\to \Theta(A(X))$. Moreover, since $i: X \to A(X)$ is an isometric
embedding it follows that $E:=A(X)$ is the desired normed space.
\end{proof}

\begin{remark} \ben
\item
This result in fact is known; (at least for group actions) it can
be derived from results of Pestov \cite{Pe-free}. In the
construction Arens-Eells space can be replaced by \emph{Free
Banach spaces}, as in above mentioned work of Pestov.
\item
Proposition \ref{l:AE} provides only a sufficient condition for
linearizability of contractive actions because we assume that the
metric space $(X,d)$ is bounded (which certainly is not a
necessary condition). The same restriction, as to our knowledge,
appears in each previous form of equivariant Arens-Eells embedding
(see e.g. \cite{Pe-free}). An elegant necessary and sufficient
condition has been recently found by Schr\"oder \cite{Sc}.
Precisely he shows that the contractive (\emph{non-expansive}, in
other terminology) $S$-action on $(X,d)$ is linearizable if and
only if all orbits $Sx$ ($x \in X$) are bounded. \een
\end{remark}

\sk
\section{$S$-Compactifications and functions} \label{s:fun}
\sk

Here we discuss how the classical Gelfand-Naimark 1-1 correspondence between Banach
subalgebras of $C(X)$ and the compactifications of $X$ can be
extended to the category of $S$-flows.

This theory is well known for topological group actions (see, for
example, J. de Vries \cite{vr-can, vr-embed}). One can easily
extend it to the case of semigroup actions (Ball and Hagler
\cite{BH}).

Separately continuous compactifications are closely related
 to the theory of right topological compactifications
 and $LMC$-functions (see Definition \ref{d:LMC}).

%

First we briefly recall some classical facts about
compactifications. A compactification of $X$ is a pair $(Y,\nu)$
where $Y$ is a compact (Hausdorff) space and $\nu$ is a continuous
map with a dense range. If $\nu$ is a topological embedding then
the compactification is said to be \emph{proper}.

Due to the Gelfand-Naimark theory there is a 1-1 correspondence
(up to the equivalence classes of compactifications) between
Banach {\it unital} (that is, the containing the constants)
subalgebras $\Acal \subset C(X)$ and the compactifications $\nu: X
\to Y$ of $X$.
 Any Banach unital $S$-subalgebra $\Acal$ of $C(X)$, induces
the {\it canonical $\Acal$-compactification} $\a_{\Acal}: X \to
X^{\Acal}$, where $X^{\Acal}$ is the Gelfand space
(or, the {\it spectrum} -- the set $MM(\Acal)$ of all
multiplicative means \cite{BJM})  of the algebra $\Acal$ (see also
Definition \ref{d:repr}.1). The map $\a_{\Acal}: X \to X^{\Acal}$
is defined by the {\it Gelfand transform}, the evaluation at $x$
multiplicative functional, that is $\a(x)(f):=f(x)$. Conversely,
every compactification $\nu: X \to Y$ is equivalent to the {\it
canonical $\Acal_{\nu}$-compactification} $\a_{\Acal_{\nu}}: X \to
X^{\Acal_{\nu}}$, where the algebra $\Acal_{\nu}$
is defined as the image $j_{\nu}(C(Y))$ of the embedding $
j_{\nu}: C(Y) \to C(X), \ \phi \mapsto \phi \circ \nu. $



\begin{remark} \label{r:domin}
If $\nu_1: X \to Y_1$ and $\nu_2: X \to Y_2$ are two
compactifications then $\nu_1$ {\it dominates} $\nu_2$, that is,
$\nu_1= q \circ \nu_2$ for some (uniquely defined) continuous map
$q: Y_2 \to Y_1$ iff $\Acal_{\nu_1} \subset \Acal_{\nu_2}$.
Moreover, if in addition, $\nu_1$ and $\nu_2$ are $S$-equivariant
maps and all $s$-translations on $X$, $Y_1$ and $Y_2$ are
continuous then $q$ is also $S$-equivariant.
\end{remark}



\begin{defin} \label{d:comp}
Let $X$ be an $S$-flow.
\ben
\item
A {\it semitopological $S$-compactification} of $X$ is a
continuous $S$-map $\alpha: X \to Y$  with a dense range
into a compact semitopological $S$-flow $Y$.
\item
Let $M \subset S$. We say that a semitopological
$S$-compactification $\a: X \to Y$ is $M$-{\it topological} if the
action $S \times Y \to Y$ is continuous at every $(m,y) \in M
\times Y$. If $M=S$ then we say {\it topological
$S$-compactification}.
\item
 A flow $(S, X)$ is said to be {\it compactifiable (semi-compactifiable)} if
there exists a {\it proper} topological (resp.: semitopological)
$S$-compactification $X \hookrightarrow Y$. A topological
semigroup $S$ is {\it compactifiable (semi-compactifiable)} if the
flow $(S,S)$, left regular action, is compactifiable (resp.:
semi-compactifiable).
\een
\end{defin}

\begin{defin} \label{d:semcomp} Let $S$ be a topological semigroup.
\ben
\item \cite{BJM}
 A {\it right topological
semigroup compactification} of $S$ is a pair $(T,\gamma)$ such
that $T$ is a compact right topological semigroup, and $\gamma$ is
a continuous homomorphism from $S$ into $T$, where $\gamma(S)$ is
dense in $T$ and the translation $\lambda_s: T \to T, \ x\mapsto
\gamma(s)x$ is continuous for every $s \in S$. It follows that the
{\it associated action} (the {\it associated flow} in \cite{La})
$$\pi_{\gamma}: S \times T \to T, \ \ (s,x) \mapsto
\gamma(s) x=\lambda_s(x)$$
is separately continuous.
Moreover, $\gamma: S \to T$ is a semigroup compactification iff
$\gamma$ is a semitopological $S$-compactification of the $S$-flow $S$
such that at the same time
$\gamma$ is a homomorphism of semigroups.
\item  A \emph{dynamical right topological
semigroup compactification} of $S$ in the sense of Ruppert
\cite{Ru} (see also \emph{monoidal compactification} of Lawson
\cite{La}) is a right topological semigroup compactification
$(T,\gamma)$ such that $\gamma$ is a topological
$S$-compactification. That is, the action $\pi_{\gamma}: S \times
T \to T$ is jointly continuous.
%
\een
\end{defin}

Evidently every semi-compactifiable flow, as a space, must be
Tychonoff.

\begin{defin} \label{d:ellis}
\ben
\item
The {\em enveloping (or Ellis) semigroup\/} $E(S,X)=E(X)$ of the
semitopological compact flow $(S,X)$ is defined as the closure in
$X^X$ (with its compact, pointwise convergence topology) of the
set $\breve{S}=\{\breve{s}: X \to X\}_{s \in S}$ considered as a
subset of $X^X$. With the operation of composition of maps this is
a right topological semigroup.
\item
The associated homomorphism $j: S \to E(X), \ s \mapsto \breve{s}$
is a right topological semigroup compactification of $S$. More
generally, for every semitopological $S$-flow $X$ and a
semitopological $S$-compactification $\a:X \to Y$ we have the
induced right topological semigroup compactification $j_{\a}: S
\to E(Y)$ such that the pair
$$
(j_{\a},\a): (S,X) \rightrightarrows (E(Y),Y)
$$
is equivariant. The associated action $\pi_j: S \times E(Y) \to
E(Y)$ is separately continuous. Furthermore, if $Y$ is a
topological $S$-flow then $\pi_j$ is jointly continuous.
\een
\end{defin}

\begin{prop} \label{l:to-ellis}
Let $S$ be a topological
semigroup.
\ben
\item
$S$ is compactifiable if and only if $S$ has a proper {\it
dynamical compactification}.
\item
$S$ is semicompactifiable if and only if it admits a proper right
topological semigroup compactification. \een
\end{prop}
\begin{proof} (2): \
 Let $\gamma: S \to T$ be a proper right topological
semigroup compactification of $S$. The associated action
$\pi_{\gamma}: S \times T \to T$ is separately continuous. Hence
$\gamma$ is a semitopological (proper) compactification of $S$.

Conversely, let $\a: S \to Y$ be a semitopological
$S$-compactification of $S$ (acting on itself by left
translations). We can pass, as in Definition \ref{d:ellis}, to the
right topological semigroup compactification $j_{\a}: S \to E(Y)$.
We can suppose without restriction of generality that $S$ is a
topological monoid (adjoining to $S$ an isolated identity $e_S$ if
necessary as in Remark \ref{r:trivial-cases}.1) and
$j_{\a}(e_S)=id_Y$. Then we have the continuous map $\hat{e}: E(Y)
\to Y, \ p \mapsto p(\a(e))$ such that $\hat{e} \circ j_{\a}=\a$.
It follows that if $\a$ is a proper compactification then $j_{\a}$
is also proper.

(1): Is similar. Observe that $\pi_j$ is jointly continuous if
$\a$ is a topological $S$-compactification.
\end{proof}

\begin{remark} \label{r:comp}
\ben
\item
For many natural monoids a separately continuous monoidal action
$\pi: S \times Y \to Y$ on arbitrary compact space $Y$ is
continuous at every $(e,y) \in \{e\} \times Y$. This happens for
instance if $S$ is a {\it Namioka space} (see \cite[Corollary
5]{La-additional} and \cite{La-points, He}). Every
\v{C}ech-complete (e.g., locally compact or complete metrizable)
space is a Namioka space.
It follows that if the monoid $S$ is a Namioka space then every
semitopological $S$-compactification $\a: X \to Y$ is
$\{e\}$-topological (or, equivalently, $H(e)$-topological, where
$H(e)$ denotes the group of all invertible elements in $S$.
\item
Recall also that by a result of Dorroh
\cite{Do} every separately continuous action of the one-parameter additive monoid $([0,\infty),+)$
on a locally compact space $X$ is jointly continuous.
\een
\end{remark}

The following fact is well known.

\begin{lemma} \label{l:r=d}
Let $G$ be a \v{C}ech-complete (e.g., locally compact or complete
metrizable) topological group. Then $\gamma: G \to T$ is a right
topological semigroup compactification of $G$ if and only if
$\gamma$ is a dynamical compactification of $G$.
\end{lemma}
\begin{proof}
In Definition \ref{d:semcomp}, (2) implies (1). The converse is
true for every topological group $S$ the underlying space of which
is \v{C}ech-complete (by Remark \ref{r:comp}.1).
\end{proof}

\begin{lemma} \label{l:repr}
Every continuous representation $(h, \alpha)$ of an $S$-space $X$
on a normed space $V$ induces the topological $S$-compactification
$$\a: X \to Y:=cl(\a(X)) \subset B^*$$
where $cl(\a(X))$ is the weak star closure of $\a(X)$ in $B^*$.
\end{lemma}
\begin{proof} Indeed,
by Lemma \ref{l:opp} the action $S \times B^* \to B^*$ is continuous. In particular,
the restricted action $S \times Y \to Y$ is continuous, too.
\end{proof}

\sk

The following definition is well known (under different names and
sometimes replacing "right" by "left") for topological group
actions \cite{vr-embed, vrexist} and for semigroups \cite{IM, BJM,
Dz, Ru, BH}.

\begin{defin} \label{d:ruc}
Let $\pi: S\times X \to X$ be a given action. A bounded function
$f \in C(X)$ is said to be {\it right uniformly continuous}
if the orbit map ${\tilde f}: S \to C(X)$ is continuous. Or,
equivalently, for every $s_0 \in S$ and $\varepsilon
>0$ there exists a neighborhood $U$ of $s_0$
such that $|f(sx)-f(s_0x)| < \varepsilon$ for every $(s,x) \in U
\times X$.
\end{defin}

For every $S$-flow $X$ denote by $RUC_S(X)$,
or, by $RUC(X)$ (where $S$ is understood) the set of all
functions on $X$ that are right uniformly continuous.
The set $RUC_S(X)$ is an $S$-invariant
Banach unital subalgebra of $C(X)$. If $X$ is a compact $S$-space
then the standard compactness arguments show that $C(X)=RUC_S(X)$.
%
%
If $X=S$ with the left regular action of $S$ on itself by left
translations, then we simply write $RUC(S)$. If $S=G$ is a
topological group, then $RUC(G)$ is the set of all usual right
uniformly continuous functions on $G$.



Let $\a_{\Acal}: X \to X^{\Acal}$ be the canonical $\Acal$-compactification of $X$.
If the Banach unital subalgebra $\Acal \subset C(X)$ is {\it $S$-invariant} (that is, the function
$(fs)(x):=f(sx)$ lies in $\Acal$ for every $s \in S$)
 then the spectrum $X^{\Acal} \subset \Acal^*$ admits the natural adjoint action
$S \times X^{\Acal} \to X^{\Acal}$ such that all translations
$\breve{s}: X^{\Acal} \to X^{\Acal}$ are continuous and
$\a_{\Acal}: X \to X^{\Acal}$ is $S$-equivariant. We get a
representation
$$(h,\a_{\Acal}): (S,X) \rightrightarrows (\Theta(\Acal)^{op}, B^*)$$
on the Banach space $\Acal$, where $h(s)(f):=fs$ and
$\a_{\Acal}(x)(f):=f(x)$. We call it the {\it canonical
$\Acal$-representation}. Note that this representation is not
necessarily \emph{continuous} because $h$ need not continuous.

\begin{prop} \label{p:RUC-comp}
Let $X$ be an $S$-flow. Assume that $\Acal$ is an $S$-invariant
unital Banach subalgebra of $C(X)$.
\ben
\item
$\a_{\Acal}: X \to X^{\Acal}$ is a topological (i.e. jointly
continuous) compactification of the $S$-flow $X$ if and only if
$\Acal \subset RUC_S(X)$.
\item
The compactification $\a_{RUC}: X \to X^{RUC}$
(for the algebra $\Acal:=RUC_S(X)$) is the maximal
topological compactification of the $S$-flow $X$.
\een
\end{prop}
\begin{proof}
(1): If $\Acal$ is a subalgebra of $RUC_S(X)$ then by Definition \ref{d:ruc} the orbit map
$\tilde{f}: S \to \Acal$ is norm continuous for every $f \in \Acal$. Therefore the canonical
representation
$$
(h,\a_{\Acal}): (S,X) \rightrightarrows (\Theta(\Acal), B^*)
$$
is continuous (because $h$ is continuous). By Lemma \ref{l:repr}
we get that the induced compactification $\a_{\Acal}: X \to
X^{\Acal}$ is a topological compactification of the $S$-flow $X$.

Conversely, if $\a_{\Acal}: X \to Y:=X^{\Acal}$ is a topological compactification then
$C(Y)=RUC_S(Y)$. This easily implies that $\Acal \subset RUC_S(X)$.

(2): Follows from (1) and Remark \ref{r:domin}.
\end{proof}

The maximal jointly continuous compactification $\a_{RUC}: S \to S^{RUC}$ defined
for the flow $(S,S)$ is the semigroup version of the so-called
"greatest ambit". Clearly, $S$ is compactifiable iff
$\a_{RUC}$ is a proper compactification. Every Hausdorff topological group $G:=S$
is compactifiable
because the algebra $RUC(G)$ separates points and closed subsets.
 It follows that the corresponding canonical representation
 (one may call it the {\it Teleman's representation})
$(h,\a_{RUC}): (G,G) \rightrightarrows (\Theta(V)^{op}, B^*)$ on
$V:=RUC(G)$ is proper and $h$ induces in fact a topological group
embedding of $G$ into $Is(V)$. The corresponding proper
compactification $\a_{RUC}: G \hookrightarrow G^{RUC}$ is the {\it
greatest ambit of $G$} (see, for example, \cite{Te, Br, vrbook2,
Pe1, UspComp}). The induced representation $(h,\a): (G,G)
\rightrightarrows (C(B^*,B^*), B^*)$ on the compact space $B^*$ is
also proper and $h$ induces an embedding of topological groups $G
\hookrightarrow Homeo(B^*)$.

Note that the maximal $S$-compactification $\a_{RUC}: X \to
X^{RUC}$ may not be an embedding even for Polish topological group
$S:=G$ and a Polish phase space $X$ (see \cite{Me1}); hence $X$ is
not $G$-compactifiable. If $S$ is discrete then $\beta_S
X=X^{RUC}$ coincides with the usual maximal compactification
$\beta X=X^{C(X)}$.

\begin{remark} \label{r:trivial-cases}
\ben
\item
Every topological semigroup $S$ canonically can be embedded into a
topological monoid $S_e:=S \sqcup \{e\}$ as a clopen subsemigroup
by adjoining to $S$ an isolated identity $e$. Furthermore, any
action $\pi: S \times X \to X$ naturally extended to the monoidal
action $\pi_e: S_e \times X \to X$. It is easy to check that
$RUC_{S_e}(X)=RUC_S(X)$. Therefore, $S$-space $X$ is
compactifiable iff $S_e$-space $X$ is compactifiable. Similarly,
$f \in RUC(S_e)$ iff $f|_S \in RUC(S)$. It follows that $S_e$ is
compactifiable iff $S$ is compactifiable.
\item
Let $Z:=X \sqcup Y$ be a disjoint sum of $S$-spaces. Then $f \in
RUC(Z)$ iff $f|_X \in RUC(X)$ and $f|_Y \in RUC(Y)$. It follows
that $Z$ is $S$-compactifiable iff $X$ and $Y$ are
$S$-compactifiable.
 \een
\end{remark}

\sk

Now we turn to the case of semitopological $S$-compactifications.

Let $ (h, \alpha): (S,X) \rightrightarrows (\Theta(V)^{op}, B^*) $
be a representation of a flow $(S,X)$ on a normed space $V$. Every
pair of vectors $(v,\psi) \in V \times V^*$ defines the function
$$
m_{v,\psi}: S \to \R, \quad s \mapsto \psi(vs)
$$
which is said to be a {\it matrix coefficient} of the given
$V$-representation.

\begin{lemma} \label{l:semi}
Let $V$ be a normed space, $X$ is an $S$-space and the pair
$$
(h,\a): (S,X) \rightrightarrows (\Theta(V)^{op}, B^*)
$$
is a representation ($h$ is not necessarily continuous).
 The following conditions are equivalent:
 \ben
\item The induced action $S \times Y \to Y$, where $Y:=cl(\a(X)) \subset B^*$, is separately
continuous (equivalently, $\a: X \to Y$ is a semitopological
$S$-compactification).
\item
The matrix coefficient
$
m_{v,\psi}: S \to \R
$
is continuous for every $v \in V$ and $\psi \in Y$.
\een
\end{lemma}
\begin{proof}
Observe that the orbit map $\tilde{\psi}: S \to Y$ (with $\psi \in
Y$) is weak star continuous if and only if the matrix coefficient
$m_{v,\psi}$ is continuous for every $v \in V$.
\end{proof}

%

This lemma naturally leads to the following definition which is
well known at least for the particular case of the left action of
$S$ on itself. It can be treated as a natural flow generalization
of the concept of {\it LMC-functions} introduced for semigroups by
Mitchell (see, for example, \cite{Mi, HM, BJMo, BJM}). However, in
general context of actions, this definition seems to be new even
for group actions.

\begin{defin} \label{d:LMC} (LMC-functions -- generalized version)
Let $X$ be an $S$-space. We say that a function $f \in C(X)$ is
{\it left multiplicatively continuous} (notation: $f \in
LMC_S(X)$, or simpler $f \in LMC(X)$) if for every $\psi \in
Y:=\beta X$ the matrix coefficient $ m_{f,\psi}: S \to \R
$ of the canonical $C(X)$-representation of $(S,X)$ is continuous.
\end{defin}

We omit a straightforward verification of the following lemma.

\begin{lemma} \label{l:LMC-pr} Let $X$ be an $S$-space.
 The set $LMC_S(X)$ is an $S$-invariant Banach subalgebra of
$C(X)$ and contains $RUC_S(X)$.
\end{lemma}

%
%
%

\begin{prop} \label{p:LMC}
Let $X$ be an $S$-space. Assume that $\Acal$ is an $S$-invariant
unital Banach subalgebra of $C(X)$ and $f \in \Acal$. \ben
\item
$f \in LMC_S(X)$ iff
for every $\psi \in X^{\Acal} \subset B^*$
the matrix coefficient
$
m_{f,\psi}: S \to \R
$
of the canonical $\Acal$-representation
is continuous.
\item
$\a_{\Acal}: X \to X^{\Acal}$ is a semitopological
compactification of the $S$-flow $X$ if and only if $\Acal \subset
LMC_S(X)$. That is, $S$-invariant unital closed subalgebras of
$LMC_S(X)$ correspond to semitopological $S$-compactifications of
$X$.
\item
The compactification $\a_{LMC}: X \to X^{LMC}$
(for $\Acal:=LMC_S(X)$) is the maximal
semitopological compactification of the $S$-flow $X$.
\item
$(S,X)$ is semicompactifiable iff $LMC_S(X)$
separates points and closed subsets in $X$.
\item (compare \cite[Ch. III, Theorem 4.5]{BJMo})
A topological semigroup $S$ is semicompactifiable iff $LMC(S)$
separates points and closed subsets in $S$ iff it admits a proper
right topological semigroup compactification.
 \een
\end{prop}
\begin{proof}
(1): The canonical $C(X)$-representation of $(S,X)$ induces the
usual maximal compactification $\beta: X \to \beta X$. Denote by
$\a_{\Acal}: X \to Y:=cl(\a_{\Acal}(X))$ the induced
compactification of the $\Acal$-representation $ (h,\a_f): (S,X)
\rightrightarrows (\Theta(\Acal), B^*). $ Then there exists a
continuous $S$-equivariant onto map $q: \beta X \to Y$ such that
$q \circ \beta = \a_{\Acal}$. It follows that the matrix
coefficient $m_{f,p}$ coincides with $m_{f,q(p)}$ for every $p \in
\beta X$.

(2): Combine Lemma \ref{l:semi} and the first assertion.

(3): Easily follows from (2).

(4): Follows from assertion (3).

(5): Use (4) and Proposition \ref{l:to-ellis}.2.
\end{proof}

%

%

Let $S$ be a topological semigroup. Then by results of
\cite[Chapter III]{BJMo} (or by the results of the present
section) we get in fact that the universal $LMC$-compactification
$u_{LMC}: S \to S^{LMC}$ (induced by the whole algebra $LMC(S)$)
of the $S$-flow $S$ is the universal right topological semigroup
compactification of $S$. Therefore our definitions and the
traditional \emph{semigroup approach} to LMC-compactifications
agree. Recall that if $G$ is a topological group that is a Namioka
space then $LMC(G)=RUC(G)$ (see \cite[Ch. III, Theorem
14.6]{BJMo}, \ Remark \ref{r:comp}.1 and Lemma \ref{l:r=d}).

\bsk
\section{$S$-compactifiability: necessary and sufficient conditions}
\sk

Let $(X,\mu)$ be a uniform space. Denote by $j_X$ or $j$ the
completion $(X,\mu) \to (\widehat{X}, \widehat{\mu})$. As usual,
$(X,\mu)$ is \emph{precompact} (or, \emph{totally bounded}) means
that the completion $(\widehat{X}, \widehat{\mu})$ is compact.
Every uniform structure $\mu$ contains the \emph{precompact
replica} of $\mu$. It is the finest precompact uniformity
$\mu_{fin} \subset \mu$. Denote by
$$i_{fin}: (X,\mu) \to (X,\mu_{fin}), \ x \mapsto x$$ the
corresponding uniform map. This map is a homeomorphism because
$top(\mu)=top(\mu_{fin})$.
 The uniformity $\mu_{fin}$
is separated and hence the corresponding completion $(X,
\mu_{fin}) \to (\widehat{X},\widehat{\mu_{fin}})=(uX,\mu_u)$ (or
simply $uX$) is a proper compactification of the topological space
$(X, top(\mu))$. The compactification $u_X=u_{(X,\mu)}: X \to uX$
is the well known \emph{Samuel compactification} (or,
\emph{universal uniform compactification}) of $(X,\mu)$.
 The corresponding algebra $\Acal_{\mu} \subset C(X)$ consists with all
$\mu$-uniformly continuous real valued bounded functions. Here we
collect some known auxiliary results.

\begin{lemma} \label{l:facts}
\begin{enumerate}
    \item For every uniform map $f: (X, \mu) \to (Y,\xi)$ the canonically associated maps
    $$f: (X, \mu_{fin}) \to (Y,\xi_{fin}),$$
    $$\widehat{f}: (\widehat{X},\widehat{\mu}) \to (\widehat{Y},\widehat{\xi})$$
    $$f^u: uX \to uY$$ are uniform.
    \item $u_X: X \to uX$ and $u_{\widehat{X}} \circ j: X \to u\widehat{X}$
    (for $u_{\widehat{X}}: \widehat{X} \to u\widehat{X}$)
    are equivalent compactifications.
    More precisely, there exists a unique homeomorphism
    $j^u: uX \to u\widehat{X}$ such that $j^u \circ u_X=u_{\widehat{X}} \circ j$.
    In particular, the natural uniform map
    $$\phi_X:=(j^u)^{-1} \circ u_{\widehat{X}}: \widehat{X} \to uX$$
    is a topological embedding.
    \item $$Unif(X,X) \to Unif(\widehat{X},\widehat{X}), \quad f \mapsto \widehat{f}$$ is a uniform embedding,
    $$Unif(X,X) \to Unif(X_{fin}, X_{fin}), \quad f \mapsto f$$ and
    $$Unif(X,X) \to Unif(uX,uX), \quad f \mapsto f^u$$ are uniform injective maps.
\end{enumerate}
\end{lemma}
\begin{proof} (1) and (3) are straightforward. For (2) observe
that the natural map $$Unif(X,\R) \to Unif(\widehat{X},\R), \ f
\to \widehat{f}$$ is a topological isomorphism of Banach algebras.
It follows that the compactifications  $u_X: X \to uX$ and
$u_{\widehat{X}} \circ j: X \to u\widehat{X}$
    are equivalent.
\end{proof}

Another direct proof of the fact that $\phi_X: \widehat{X} \to uX$
is a uniform embedding can be found in \cite{FS}.

\begin{defin} \label{d:equiun}
Let $\mu$ be a uniformity on $X$ and $\pi: S\times X \to X$ be a semigroup action.
We call this action:
\ben
\item
{\it $\mu$-saturated} if every $s$-translation $\breve{s}: X \to
X$ is $\mu$-uniform (thus the corresponding homomorphism $h_{\pi}:
S \to Unif(X,X), \ s \mapsto \breve{s}$ is well defined).
\item
 {\it $\mu$-bounded at $s_0$} if for every $\varepsilon \in \mu$
 there exists a neighborhood $U(s_0)$ such that $(s_0x,sx) \in
\varepsilon$ for each $x\in X$ and $s\in U$. If this condition
holds for every $s_0 \in S$ then we simply say:
\emph{$\mu$-bounded.}
\item (see \cite{me-gr-com}) \emph{$\mu$-equiuniform} if $\mu$ is saturated and bounded.
It is equivalent to say that the corresponding homomorphism
$h_{\pi}: S \to Unif(X,X)$ is continuous.
\item \emph{$(\xi,\mu)$-equiuniform} if
$\xi$ is a compatible uniformity on $S$ such that the left actions
$\nu: S \times S \to S$ and $\pi: S\times X \to X$ are saturated
(with respect to $\xi$ and $\mu$ respectively) and the associated
homomorphisms $h_{\pi}: S \to Unif(X,X)$, \ $h_{\nu}: S \to
Unif(S,S)$ are uniform maps. \een

Sometimes we say also that the uniformity $\mu$ is saturated,
bounded and equiuniform, respectively.
\end{defin}

For group actions bounded uniformities appear in \cite{vr-embed}
and in \cite{Ca} (see also "uniform action" in the sense of
\cite{Ak}). We collect here some simple examples.

\begin{examples} \label{e:easy}
\ben
\item
Every $\mu$-equiuniform action is continuous.
 \item
 Every compact $S$-space $X$ is equiuniform (with respect to the unique compatible uniformity on $X$).
\item
For every uniform space $(X,\mu)$ and every subsemigroup $S
\subset Unif(X,X)$ endowed with the subspace uniformity $\xi$
inherited from $Unif(X,X)$ the natural action $S \times X \to X$
(see Example \ref{ex:monoids}.1) is $(\xi,\mu)$-equiuniform.
%
%

\item
For every $(\xi,\mu)$-equiuniform action $S \times X \to X$ the
left action $S \times S \to S$ is $(\xi,\xi)$-equiuniform.

\item Let $S$ be a semigroup with a right invariant
uniformity $\xi$ on $S$ such that all left translations are
uniformly continuous. Then the left action $S \times S \to S$ is
$(\xi,\xi)$-equiuniform.



\een
\end{examples}


We need some notation. Let $S \times X \to X$ be a semigroup
action. For every element $s \in S$ and a subset $A \subset X$
define $s^{-1}A:=\{x \in X: sx \in A \}$. Let $\mu$ be a
uniformity on $X$ and $\varepsilon \in \mu$. Then $\eps$ is a
subset of $X \times X$. For every $s \in S \cup \{id_X\}$ we can
define similarly the following set
$$s^{-1}\eps:=\{(x,y) \in X\times X: (sx,sy) \in \varepsilon\}$$
where $id_X^{-1} \eps=\eps$.

\begin{lemma} \label{l:noteasy}
Let $\mu$ be a uniformity on $X$ such that the
semigroup action of a topological semigroup $S$ on $(X, top(\mu))$
is continuous. \ben
\item
The family $\{s^{-1}\varepsilon: s\in S \cup \{id_X\}, \varepsilon \in \mu\}$
is a subbase of a saturated uniformity $\mu^S \supseteq \mu$
generating the same topology (that is, $top(\mu)=top(\mu^S)$).
\item
If the action is $\mu$-bounded
then it is also $\mu^S$-bounded (hence, $\mu^S$-equiuniform).
\item
If the action is $\mu$-bounded ($\mu$-saturated, $\mu$-equiuniform, or $(\xi,\mu)$-equiuniform)
then it is also $\mu_{fin}$-bounded ($\mu_{fin}$-saturated, $\mu_{fin}$-equiuniform,
or $(\xi,\mu_{fin})$-equiuniform
respectively).
\een
\end{lemma}
\begin{proof}
The proofs of (1) and (2) are trivial.

(3): The boundedness of $\mu_{fin}$ is clear because $\mu_{fin}
\subset \mu$. In order to show that the action is
$\mu_{fin}$-saturated we have to check that $\tilde{s}: (X,
\mu_{fin}) \to (X, \mu_{fin})$ is uniform for every $s \in S$. Let
$\eps \in \mu_{fin}$. Since $s(s^{-1}\eps) \subset \eps$ we have
only to show that $s^{-1}\eps \in \mu_{fin}$.

Pick a symmetric entourage $\delta \in \mu_{fin}$ such that
$\delta \circ \delta \subset \eps$. Since $\delta \in \mu_{fin}$
there exists a finite
 subset $\{y_1,y_2, \cdots, y_n\}$ in $X$ which is $\delta$-dense in $X$
 (that is, $\cup_{i=1}^{n} \delta(y_i) = X$,
 where $\delta(y):=\{x \in X: \ (x,y) \in
 \delta\}$). Passing to a subfamily if necessary we can suppose in addition
that $\delta(y_i) \cap sX \neq \emptyset$ for every $i \in
\{1,2,\cdots,n\}$. Choose $z_i \in X$ such that $sz_i \in
\delta(y_i)$ for each $i$. Then $\{z_1,z_2, \cdots, z_n\}$ is a
 finite $s^{-1}\eps$-dense subset in $X$.
Indeed, for every $x_0 \in X$ there exists $i_0$ such that
$(y_{i_0},sx_0) \in \delta$. Since $(sz_{i_0},y_{i_0}) \in \delta$
we get $(sz_{i_0}, sx_0) \in \delta \circ \delta \subset \eps$.
Thus, $(z_{i_0}, x_0) \in s^{-1}\eps$.

%

Checking that the action is $(\xi,\mu_{fin})$-equiuniform
(provided that it is $(\xi,\mu_{fin})$-equiuniform) observe that
the map $Unif(X,X) \to Unif(X_{fin}, X_{fin}), \ f \mapsto f$ is
uniform. This implies that the homomorphism $(S,\xi) \to
Unif(X_{fin}, X_{fin})$ is also uniform.
\end{proof}

\begin{lemma} \label{l:completion1}
\ben
\item Let $\mu$ be a saturated uniformity on $X$
with respect to the action $S \times X \to X$.
Let $Y$ be an $S$-invariant dense subset of $X$
such that the induced action
$S \times Y \to Y$ is $\mu|_Y$-bounded.
Then the given action $S \times X \to X$ is $\mu$-equiuniform and
continuous.
\item
Let $\pi: S\times X \to X$ be a continuous $\mu$-equiuniform action.
Then the induced action on the completion
${\widehat \pi}: S\times \widehat X \to {\widehat X}$ is well-defined, ${\widehat \mu}$-equiuniform
(and continuous).
\een
\end{lemma}
\begin{proof}
(1) Let $s_0 \in S$ and $\eps \in \mu$. There exists
an element $\eps_1 \in \mu$ such that $\eps_1 \subset \eps$ and
$\eps_1$ is a closed subset of $X \times X$. Choose a neighborhood
$U(s_0)$ such that $(s_0y,sy) \in \eps_1$ for every $s \in U$ and
$y \in Y$. Then $(s_0x,sx) \in \eps$ for every $s \in U$ and $x
\in X$. Thus the given (saturated) action is $\mu$-bounded. The
action is continuous by Example \ref{e:easy}.1.

(2) Easily follows from (1).
\end{proof}


\begin{lemma} \label{l:extensions}
Let $X$ and $P$ be Hausdorff spaces. Assume that:
\begin{itemize}
    \item [(i)] $S$ is a dense subset of $P$.
    \item [(ii)] $S$ is a semigroup w.r.t. the operation $w_S: S \times S \to S$.
    \item [(iii)] $\vartheta: S \times P \to P$ is a semigroup action with continuous translations.
    \item [(iv)] $m: P \times P \to P$ is a right continuous mapping which extends $w_S$ and $\vartheta$.
    \item [(v)] $\pi_S: S \times X \to X$ is a semigroup action with continuous translations.
    \item [(vi)]$\pi_P: P \times X \to X$ is a right continuous mapping which extends $\pi_S$.
\end{itemize}

Then we have:
\begin{enumerate}
    \item $(P,m)$ is a right topological semigroup.
    \item $\pi_P: P \times X \to X$ is a semigroup action.
    \item If $X$ is regular and $\pi_S$ is continuous
    at $(s_0,x_0)$ with some $(s_0, x_0) \in S \times X$
    then $\pi_P$ remains continuous at $(s_0,x_0)$.
\end{enumerate}
\end{lemma}
\begin{proof} First of all we check the associativity
$$(p_1p_2)x=p_1(p_2x)$$
for every given triple $(p_1, p_2,x) \in  P \times P \times X$,
where $(p_1p_2)x:=\pi_P (m(p_1,p_2),x)$ and
$p_1(p_2x):=\pi_P(p_1,\pi_P(p_2,x))$.

Choose nets $a_i$ and $b_j$ in $P$ such that $a_i, b_j \in S$ and
$\lim_i a_i =p_1, \lim_j b_j=p_2$. Then by our assumptions we have
$(p_1p_2)x=\lim_i (a_ip_2)x=(\lim_i \lim_j (a_ib_j))x = \lim_i
\lim_j (a_i(b_jx))=\lim_i a_i(\lim_j (b_jx))=\lim_i
(a_i(p_2x))=p_1(p_2x).$

Apply this formula to the particular case of $X:=P$. Then we get
that $(p_1p_2)p_3)=p_1(p_2p_3)$ for all triples $(p_1,p_2,p_3) \in
P^3$. This proves (1). Moreover, now the general formula means
that $\pi_P$ is a semigroup action.

For $(3)$ use Lemma \ref{l:contpoints}.
\end{proof}

\begin{prop} \label{p:semcompletion}
Let $\xi$ be a compatible uniformity on a topological semigroup
$S$  such that the left action $\nu: S \times S \to S$ is $(\xi,
\xi)$-equiuniform. Identify $S$ with its image under the
completion map $j: S \to \widehat{S}$. Then there exists a map $m:
\widehat{S} \times \widehat{S} \to \widehat{S}$ such that
$(\widehat{S},m)$ is a topological semigroup, $S$ is a
subsemigroup of $\widehat{S}$ and the left action $m$ is
$(\widehat{\xi},\widehat{\xi})$-equiuniform.
\end{prop}
\begin{proof}
The natural homomorphism $h_{\nu}: (S,\xi) \to Unif(S,S), \ s
\mapsto \lambda_s$ is uniform. Consider the uniform embedding
$$Unif(S,S) \to Unif(\widehat{S},\widehat{S}), \quad f \mapsto \widehat{f}.$$

Denote by $h$ the corresponding uniform composition $h: S \to
Unif(\widehat{S},\widehat{S})$. Since the uniform space
$Unif(\widehat{S},\widehat{S})$ is complete there exists a unique
uniform extension $\widehat{h}: \widehat{S} \to
Unif(\widehat{S},\widehat{S})$ of $h$. Then the evaluation map $m:
\widehat{S} \times \widehat{S} \to \widehat{S}, \
m(t,p)=\widehat{h}(t)(p)$ is jointly continuous and extends the
original multiplication $\nu$ on $S$. On the other hand by Lemma
\ref{l:completion1} we get that there exists a uniquely determined
continuous semigroup action $\vartheta: S \times \widehat{S} \to
\widehat{S}$ which also extends $\nu$. It follows that $m$ extends
$\vartheta$. By Lemma \ref{l:extensions} (for the setting
$P:=\widehat{S}, X:=\widehat{S}$) we obtain that $(\widehat{S},m)$
is a semigroup and $S$ is its subsemigroup. Furthermore,
$\widehat{S}$ is a topological semigroup because $m$ is
continuous. Since $h_{\nu}$ is a uniform homomorphism and
$\widehat{h}|_S=h_{\nu}$ it follows that the uniform map
$\widehat{h}$ also is a {\it homomorphism} of semigroups. This
means that the left action $m$ is
$(\widehat{\xi},\widehat{\xi})$-equiuniform.
\end{proof}

\begin{prop} \label{p:completion2}
Let $\pi: S \times X \to X$ be a $(\xi,\mu)$-equiuniform action.
Then there exist (uniquely determined)
continuous semigroup actions: \bit
\item [(i)] 
$\widehat{\pi}: \widehat{S} \times \widehat{X} \to \widehat{X}$
which is $(\widehat{\xi}, \widehat{\mu})$-equiuniform
and naturally extends $\pi$;
\item [(ii)]
$\pi: S \times X_{fin} \to X_{fin}$ which is $(\xi,\mu_{fin})$-equiuniform;
\item [(iii)]
$\widehat{\pi}_u: \widehat{S} \times uX \to uX$ which is $(\widehat{\xi}, {\mu_u})$-equiuniform and
naturally extends $\widehat{\pi}$.
\eit
\end{prop}
\begin{proof} (i) By Proposition \ref{p:semcompletion}
we know that the left action is
$(\widehat{\xi},\widehat{\xi})$-equiuniform on the topological
semigroup $\widehat{S}$. Since $Unif(X,X) \to
Unif(\widehat{X},\widehat{X}), f \mapsto \widehat{f}$ is a uniform
embedding and $Unif(\widehat{X},\widehat{X})$ is complete there
exists a (unique) uniform map $\widehat{h}: \widehat{S} \to
Unif(\widehat{X},\widehat{X})$ which extends the homomorphism
$h=h_{\pi}: S \to Unif(X,X)$. In fact $\widehat{h}$ is a
homomorphism because $h$ and $\widehat{h}$ agree on a a dense
subsemigroup $S$ of $\widehat{S}$ and
$Unif(\widehat{X},\widehat{X})$, $\widehat{S}$ are Hausdorff
topological semigroups. This proves that the action
$\widehat{\pi}$ is $(\widehat{\xi}, \widehat{\mu})$-equiuniform.
The action $\widehat{\pi}$ extends the original action $\pi$
because $\widehat{h}$ extends $h$.


$(ii)$ Is clear by Lemma \ref{l:noteasy}.3.

$(iii)$ Combine $(i)$ and $(ii)$ taking into account that $uX$ is
the completion of $\mu_{fin}$.

The continuity of these actions are trivial by Example
\ref{e:easy}.1.


\end{proof}

\begin{prop} \label{p:samuel}
\ben
\item
If the semigroup action $\pi: S\times X \to X$ is $\mu$-equiuniform then
the induced action $\pi_u: S\times uX \to uX$
on the Samuel compactification
$uX:=u(X,\mu)$ is a proper
$S$-compactification of $X$.
\item
$(S,X)$ is 
compactifiable iff the action on
$X$ is $\mu$-bounded
with respect to some compatible uniformity $\mu$.
\een
\end{prop}
\begin{proof} (1) The action is $\mu$-equiuniform means that the
homomorphism $h_{\pi}: S \to Unif(X,X)$ is continuous. It suffices
to prove our assertion for the action of $h_{\pi}(S) \times X \to
X$. Hence we can suppose that in fact $S$ is the semigroup
$h_{\pi}(S)$. Now the action is $(\xi,\mu)$-equiuniform where
$\xi$ is the uniformity induced on $h_{\pi}(S)$ from $Unif(X,X)$.
Using Proposition \ref{p:completion2}(iii) we get a continuous
action $\widehat{\pi}_u: \widehat{S} \times uX \to uX$ which is
$(\widehat{\xi}, {\mu_u})$-equiuniform and naturally extends
$\widehat{\pi}$. Then its restriction $\pi_u: S\times uX \to uX$
is continuous, too. Hence $u: X \to uX$ is a (proper)
$S$-compactification of $X$.

(2)
Assume that $X$ is $\mu$-bounded. Then by Lemma \ref{l:noteasy}
the action is $\mu^S$-equiuniform (which is a compatible
uniformity). Now by the first assertion $X$ is $S$-compactifiable.
For the converse use Example \ref{e:easy}.2.
\end{proof}

\begin{corol} \label{c:equiun}
There exists a 1-1 correspondence between proper topological
$S$-compactifications of $X$ and precompact compatible
equiuniformities on $X$.
\end{corol}

Note that Corollary \ref{c:equiun} is well known for group actions \cite{Br,me-gr-com}.

\begin{thm} \label{t:complet}
Let $\pi: S\times X \to X$ be a $(\xi,\pi)$-equiuniform semigroup
action.
%
%
Then
 \bit
\item [(a)]
$u: S \to uS$ is a proper right topological semigroup compactification of $S$.
\item [(b)]
There exists a right continuous semigroup action $\pi^u_u: uS \times uX \to uX$ which extends
the action $\widehat{\pi}_u: \widehat{S} \times uX \to uX$ (hence also
$\widehat{\pi}: \widehat{S} \times \widehat{X} \to \widehat{X}$)
and is continuous at every $(p,z) \in \widehat{S} \times
uX$.

\eit
\end{thm}
\begin{proof} \ By Proposition \ref{p:completion2}$(iii)$ there exists a continuous action
$\widehat{\pi}_u: \widehat{S} \times uX \to uX$ which extends
$\widehat{\pi}$ and is $(\widehat{\xi}, {\mu_u})$-equiuniform.
Then, in particular, every orbit map $\tilde{z}: \widehat{S} \to
uX, \ t \mapsto tz$ is uniform. By the universality of Samuel
compactifications there exists a uniquely defined continuous
extension $u\widehat{S} \to uX$ of $\tilde{z}$. The
compactifications $S \to uS$ and $S \to u\widehat{S}$ are
naturally equivalent (Lemma \ref{l:facts}.2). Hence we have a
continuous function $\tilde{z}_u: uS \to uX$ which extends the map
$\tilde{z}: \widehat{S} \to uX$, where $\widehat{S}$ is treated as
a topological subspace of $uS$.

Now we define $\pi^u_u: uS \times uX \to uX$ by $\pi_u^u
(p,z):=\tilde{z}_u(p)$ for every $p \in uS$ and $z \in uX$.
Clearly, $\pi^u_u$ is right continuous and
$\pi_u^u(t,z)=\widehat{\pi}_u (t,z)$ for every $t \in
\widehat{S}$. On the other hand again by Proposition
\ref{p:completion2}$(iii)$ (for $X:=S$) we have the continuous
action $\widehat{S} \times uS \to uS$ which extends the
multiplication $\widehat{m}: \widehat{S} \times \widehat{S} \to
\widehat{S}$ (via the natural dense embedding
$\widehat{S}=\phi_S(\widehat{S}) \hookrightarrow uS$). We can
apply Lemma \ref{l:extensions} (for the dense subset
$\widehat{S}=\phi_S(\widehat{S})$ of $uS$ and natural maps
$\widehat{\pi}_u$ and ${\pi}^u_u$). It follows that $uS$ is a
right topological semigroup with the subsemigroup $\widehat{S}$
and $\pi^u_u: uS \times uX \to uX$ is a right continuous semigroup
action extending $\widehat{\pi}_u$. By Lemma \ref{l:extensions}.3
we get that $\pi^u_u$ is jointly continuous at every $(p,z) \in
\widehat{S} \times uX$.
\end{proof}

\begin{corol} \label{c:str}
Let $S$ be a semigroup with a right invariant uniformity
$\xi$ on $S$ such that all left translations are
uniformly continuous.
\ben
\item \emph{(Kocak and Strauss \cite{KS})}
$S \to uS$ is a right topological semigroup compactification of $S$.
\item \emph{(Ferri and Strauss \cite{FS})}
The multiplication $uS \times uS \to uS$ is jointly continuous at every $(p,z) \in \widehat{S} \times uS$.
\een
\end{corol}
\begin{proof}
By Example \ref{e:easy}.5 the left action $S \times S \to S$ is
$(\xi,\xi)$-equiuniform. Now apply Theorem \ref{t:complet}.
\end{proof}

\sk

Now we give a compactifiability criteria for semigroup actions.

\begin{thm} \label{t:flows}
For every $S$-space $X$ the following conditions are equivalent:
\begin{enumerate}
\item $X$ is $S$-compactifiable.
\item
$RUC_S(X)$ separates points from closed subsets.
\item There exists a Banach space $V$
and a proper continuous representation
$$
(h,\a): (S,X) \rightrightarrows (\Theta(V)^{op}, B^*).
$$
\item There exists a compact space $Y$ and a proper representation
$$(h,\a): (S,X) \rightrightarrows (C(Y,Y), Y).$$
\item There exists a uniform space $(Y, \mu)$ and a proper representation
$$(h,\a): (S,X) \rightrightarrows (Unif(Y,Y), Y).$$
\end{enumerate}
\end{thm}
\begin{proof} (1) $\Rightarrow$ (2): Let $\nu: X \hookrightarrow Y$ be a proper $S$-compactification.
Then $C(Y)=RUC_S(Y)$. Now use the obvious hereditarity property of
right uniformly continuous functions. That is the fact that $f
\circ \nu  \in RUC_S(X)$ for every $f \in RUC_S(Y)$.

(2) $\Rightarrow$ (3): Consider the canonical $
V:=RUC_S(X)$-representation of $(S,X)$ on $V$ and apply
Proposition \ref{p:RUC-comp}.

(3) $\Rightarrow$ (4): Apply Lemma \ref{l:opp} to $V:=RUC_S(X)$.

(4) $\Rightarrow$ (5): For a compact space $K$ (and its uniquie
compatible uniformity) the uniform spaces $Unif(K,K)$ and $C(K,K)$
are the same.

(5) $\Rightarrow$ (1): By Example \ref{e:easy}.3 there exists a
compatible uniformity $\mu$ on $X$ such that the action is
$\mu$-equiuniform. Then the corresponding Samuel compactification
of $(X,\mu)$ is an $S$-compactification by virtue of Proposition
\ref{p:samuel}.1.
\end{proof}


The following theorem shows that a topological semigroup $S$ is compactifiable iff
$S$ "lives in natural monoids".

\begin{thm} \label{t:sem}
Let $S$ be a topological semigroup. The following are equivalent:
\ben
\item $S$ is compactifiable;
\item $RUC(S)$ determines the topology of $S$.
\item The monoid $S_e$ (from Remark
\ref{r:trivial-cases}.1) is compactifiable;
\item $S$ has a proper {\it dynamical
compactification}.
\item $S^{op}$ (the opposite semigroup of $S$)
is a topological subsemigroup of $\Theta(V)$ for some normed
(equivalently, {\it Banach}) space $V$;
\item $S^{op}$ is a
topological subsemigroup of $\Theta(M,d)$ for some
metric space $(M,d)$;
\item $S$ is a topological subsemigroup of $C(Y,Y)$
for some compact space $Y$;
\item $S$ is a topological subsemigroup of $Unif(Y,Y)$
for some uniform space $(Y, \mu)$.
\item
There exists a compatible right invariant uniformity $\mu$ on $S$.
\item
There exists a compatible uniformity $\mu$ on $S$ such that the
right action of $S$ on $(S, \mu)$ is equicontinuous.
\item The topology of $S$ can be generated by
a family $\{d_i\}_{i \in I}$ of right contractive pseudometrics on
S.
\een
If $S$ is a monoid then we can ensure in the assertions
\emph{(5), (6), (7)} and \emph{(8)} that $S$ is a topological
submonoid of the corresponding topological monoid.
\end{thm}
\begin{proof}
(1) $\Leftrightarrow$ (2): Follows from Proposition
\ref{p:RUC-comp}.

(1) $\Leftrightarrow$ (3): See Remark \ref{r:trivial-cases}.1.

(2) $\Leftrightarrow$ (4): $RUC(S)$ determines the topology iff
the universal dynamical compactification $u_{RUC}: S \to S^{RUC}$
is {\it proper}.

(1) $\Rightarrow$ (5): First of all observe that by Remark
\ref{r:simple}.2, "normed" and "Banach" cases of (5) are
equivalent.

 By our assumption $S$
is $S$-compactifiable. Theorem \ref{t:flows} implies that there
exists a proper continuous representation
$$
(h,\a): (S,S) \rightrightarrows (\Theta(V)^{op}, B^*).
$$
Where $V:=RUC(S)$. By (1) $\Leftrightarrow$ (3) we can assume that
$S$ is a monoid. Since $\a: S \to B^*$ is an $S$-embedding and the
pair $(h, \a)$ is equivariant it follows that the homomorphism $h:
S \to \Theta(V)^{op}$ is a topological embedding, too.


(5) $\Rightarrow$ (6): $\Theta(V)$ is embedded into $\Theta(M,d)$
where $M:=V$ and $d(x,y):=||x-y||$.


(5) $\Rightarrow$ (7): Immediate by Lemma \ref{l:opp}.

(7) $\Rightarrow$ (8): Trivial.


(8) $\Rightarrow$ (9): Follows by Example \ref{e:unif-eq}.2.

(9) $\Rightarrow$ (10): Trivial by Definition \ref{d:equic}.


(10) $\Rightarrow$ (11): If a family $\{d_i\}$ of bounded
pseudometrics generates an equicontinuous uniformity $\mu$ then
the family $\{d^*_i\}$ of right contractive pseudometrics
$$d^*_i(x,y):=max \{sup_{s \in S} d_i(xs,ys),
d(x,y) \}$$ generates a
uniformity $\mu^*$ which is topologically equivalent to $\mu$.

(11) $\Rightarrow$ (1): Let $\mu$ be the uniformity generated by
the given family of pseudometrics on $S$. Since the pseudometrics
are right contractive it follows that the action of $S$ on $S$ is
$\mu$-bounded. Now Proposition \ref{p:samuel}.2 implies that $S$
is a compactifiable $S$-flow.

(6) $\Rightarrow$ (11): Denote by $(S, \ast)$ the opposite
semigroup $\Theta(M,d)^{op}$ of $\Theta(M,d)$. The family of
pseudometrics $\{\rho_m\}_{m \in M}$ generates the topology of $S$
where
$$\rho_m (s_1,s_2):= d(s_1m, s_2m).$$
Now observe that each $\rho_m$ is right contractive on the
topological semigroup $S$. Indeed, for every triple $t,s_1,s_2 \in
S$ we have
$$\rho_m(s_1 \ast t,s_2 \ast t)=\rho_m (ts_1,ts_2)=d(ts_1m,ts_2m)
\leq d(s_1m, s_2m)=\rho_m(s_1,s_2).$$

\sk

Finally, note that if $S$ is a monoid then by the proof of (2)
$\Rightarrow$ (5) the homomorphism $h: S \to \Theta(V)^{op}$ is a
topological embedding of monoids.
\end{proof}

\begin{corol} \label{c:mon}
Each of the following semigroups is compactifiable:
\ben
\item
$\Theta(X,d)^{op}$ for every
metric space $(X,d)$. In
particular, $\Theta(V)^{op}$ (endowed with the strong operator
topology) for every normed space $V$.
\item
$Unif(Y,Y)$ for every uniform space $(Y, \mu)$.
\item
$C(Y,Y)$ for every compact space $Y$.
\item
$(B_V, \tau_u)$ endowed with the uniform topology for every normed
algebra $V$ (e.g., for the algebra $V:=L(E)$
for arbitrary normed space $E$).
\item
Let $G$ be a topological group and $\mathcal R$ its right
uniformity. Then the completion $S:=(\widehat{G},
\widehat{\mathcal R})$ is a topological semigroup and this
semigroup is compactifiable. \een
\end{corol}
\begin{proof} All assertions easily follow from Theorem \ref{t:sem}. For (4) observe that
the original metric of the original norm on $B_V$ is right (and
also left) contractive $||xs-ys|| \leq ||x-y|| \cdot ||s|| \leq
||x-y||$ for every $x, y, s \in B_V$.

(5): $G \times G \to G$ is $\mathcal R(G)$-equiuniform. Apply now
Propositions \ref{p:semcompletion} and \ref{p:completion2}.
\end{proof}

It is well known that $(\widehat{G}, \widehat{\mathcal R})$ is a
topological semigroup (see for example \cite[Proposition
10.12(a)]{RD}) containing $G$ as a subsemigroup. For several
important semigroups of the form $S:=(\widehat{G},
\widehat{\mathcal R})$ see Pestov \cite{Pe-nbook}.

\begin{remark} \label{r:results}
\ben
\item
Kocak and Strauss proved in  \cite[Theorem 14]{KS} that if a
topological semigroup $S$ admits a right invariant left saturated
uniformity then $S$ is compactifiable. One can remove
``saturated'' as Theorem \ref{t:sem} shows. Furthermore by
assertion (9) the existence of right invariant uniformity is also
a necessary condition.
\item
As we already have seen $\Theta(E)^{op}$ is compactifiable for
every normed space $E$. It is not true for $\Theta(E)$, in
general, as we will see later in Examples \ref{ex:count}. So we
cannot substitute $\Theta(E)^{op}$ by $\Theta(E)$ in Theorem
\ref{t:sem}. However, we can repair this situation for involutive
subsemigroups $S$ of $\Theta(E)$ (see Corollary \ref{c:inv}).
\item
We cannot change $B_V$ by $V$ in Corollary \ref{c:mon}.4 as the
example of the multiplicative semigroup $V:=\R$ shows (see
Examples \ref{ex:count}.2).
\item Our results suggest a semigroup version of the right
uniformities $\mathcal R(S)$. For a compactifiable topological
semigroup one can define $\mathcal R(S)$ as the finest right
invariant compatible uniformity on $S$. Then Corollary
\ref{c:mon}.5 admits a natural semigroup generalization for the
completion of $(S, \mathcal R(S))$.
 \een
\end{remark}

\begin{thm} \label{para}
Let $G$ be a paratopological group. Then $G$ is compactifiable iff
$G$ is a topological group.
\end{thm}
\begin{proof}
If $G$ is compactifiable then by Theorem \ref{t:sem} we have an
embedding $h: G \to C(K,K)$ of topological monoids. Then $h(G)
\subset Homeo(K)$, where $Homeo(K)$ is a topological group. The
converse is clear by the Teleman's representation.
\end{proof}



\sk

Recall that a semigroup $S$ is said to be an {\it inverse
semigroup} if for every $s \in S$ there exists a unique $s^* \in
S$ such that $ss^*s=s$ and $s^*ss^*=s^*$. {\it Topological inverse
semigroup} will mean that the multiplication is continuous and in
addition the map $S \to S, \ s\mapsto s^*$ is continuous.

By an {\it involution} on a semigroup $S$ we mean a map $i: S\to
S$ such that $i(i(s))=s$ and $i(s_1s_2)=i(s_2)i(s_1)$. If $S$
admits a continuous involution then we say that $S$ is {\it
topologically involutive}. 
For example, $S$
is involutive if $S$ is a topological inverse semigroup;. This
happens in particular if either $S$ is a commutative topological
semigroup or a topological group.

\begin{prop} \label{c:inv}
Let $S$ be a topological subsemigroup
of $\Theta(E)$ for a normed space $E$. Suppose that $S$ is
topologically involutive. Then $S$ is compactifiable.
\end{prop}
\begin{proof}
Use Corollary \ref{c:mon}.1
\end{proof}


\bsk
\section{A universal compactifiable semigroup}
\sk

Denote by $U$ the topological monoid $C(I^{\omega}, I^{\omega})$,
where $I:=[0,1]$ is the closed interval. Theorem \ref{t:sem}
implies that $U$ is compactifiable. It contains the subgroup
$Homeo(I^{\omega})$ of all selfhomeomorphisms of the Hilbert cube
$I^{\omega}$. Recall that $Homeo(I^{\omega})$ is a universal
second countable topological group (see Uspenskij \cite{Us1}).
Moreover, by \cite{Me2} the group action $Homeo(I^{\omega}) \times
I^{\omega} \to I^{\omega}$ is universal for all second countable
compactifiable $G$-flows $X$ with a second countable acting group
$G$. We can now give a natural generalization for semigroups and
semigroup actions.

\begin{thm} \label{t:unflows}
Let $S$ be a compactifiable second countable semigroup. Then every
compactifiable second countable $S$-flow $X$ is a part of the flow
$(U, I^{\omega})$. That is, there exists a representation $(h,\a):
(S,X) \rightrightarrows (U, I^{\omega})$ such that $h: S
\hookrightarrow U$ is an embedding of topological semigroups and
$\a: X \hookrightarrow I^{\omega}$ is a topological embedding.
\end{thm}
\begin{proof} By Remark \ref{r:trivial-cases}.1
we can assume that $S$ is a monoid with the identity $e$ and $S
\times X \to X$ is a monoidal action.

Furthermore, we can suppose in addition that the action is {\it
topologically exact}. This means (see \cite{Me2}) that: (a) $sx
=x$ for all $x \in X$ implies that $s=e$; (b) there exists no
strictly weaker topology on $S$ which makes the action on $X$
continuous. Indeed, we can pass, if necessary, to the following
new (but still $S$-compactifiable by Remark
\ref{r:trivial-cases}.2) second countable phase space $X':=X
\sqcup S$, a disjoint sum of the $S$-flows $X$ and $S$, where the
monoid $S$ acts on itself by left multiplications. Thus, by our
assumption $X$ is a compactifiable $S$-flow with the topologically
exact action. The algebra $RUC(X)$ separates points and closed
subsets of $X$. Since $X$ is second countable we can choose a
separable closed subalgebra $\Acal$ of $RUC(X)$ having the same
property. Moreover since $S$ is also second countable we can
assume that $\Acal$ is even $S$-invariant. Indeed if $T \subset
\Acal$ and $S_1 \subset S$ are countable dense subsets then $TS_1$
is a countable dense subset in the $S$-invariant closed subalgebra
$\Acal' \supseteq RUC(X)$ topologically generated by $\Acal S$.

Now consider the corresponding representation
$$(h,\a): (S,X) \rightrightarrows (\Theta(\Acal)^{op}, B^*)$$
of the flow $(S,X)$ on the Banach space $\Acal$. Now, as in
\cite{Us1}, we use the fact that the unit ball $B^*$ being a
convex compact subset of a separable Frechet space $(\Acal,
weak^*)$ is homeomorphic by Keller's theorem (see for example
\cite{BP}) to the Hilbert cube $I^{\omega}$. By our assumption
$\Acal$ separates points from closed subsets in $X$. Therefore the
map $\a: X \hookrightarrow B^*$ is a topological embedding.
Moreover, since the action of $S$ on $X$ is topologically exact
and the pair $(h,\a)$ is equivariant it follows that the
homomorphism $h: S \to \Theta(\Acal)^{op}$ is in fact an embedding
of topological monoids. Observe that
$$(\gamma, id): (\Theta(\Acal)^{op}, B^*) \rightrightarrows (C(B^*,B^*),B^*)$$
is an equivariant pair with the embedding $\gamma$ of topological
monoids (see Lemma \ref{l:opp}). Now substituting $B^*$ by the
Hilbert cube $I^{\omega}$ we complete the proof.
\end{proof}

As a corollary we get

\begin{thm} \label{usem}
\emph{(Semigroup version of Uspenskij's theorem)} The monoid
$U:=C(I^{\omega}, I^{\omega})$ is universal in the class of all
second countable compactifiable semigroups.
\end{thm}



\bsk

\section{Some examples}

\sk


Recall that if $G$ is a Hausdorff (Tychonoff) topological group
then a Tychonoff $G$-flow $X$ is compactifiable in each of the
following cases:

(a) $G$ is locally compact \cite{vrexist};

(b)  $X$ is locally compact \cite{vr-can};

(c) $X$ admits a $G$-invariant metric \cite{lud-vr};

(d) $X$ is a normed space and each $g$-translation $X \to X$ is
linear \cite{me-singap};

(e) $G$ is second category, $(X,d)$ is a metric $G$-space and each
$\breve{g}: X \to X$ is $d$-uniformly continuous \cite{me-singap}.

\sk

Examples below show that for the case of monoidal actions
 analogous results do not remain true, in general.

Answering de Vries' "compactification problem" negatively in
\cite{Me1} we construct a noncompactifiable Polish $G$-space $X$
with a Polish acting group $G$. Moreover by \cite{MeSc1} for every
Polish group $G$ which is not locally compact there exists a
suitable noncompactifiable Polish $G$-space. We can use this fact
below (see Example \ref{ex:count}.10) providing many
non-semi-compactifiable Polish topological semigroups. We refer to
\cite{me-singap,Me-opit2} for more information about
compactifications of group actions.

\sk

\begin{lemma} \label{key}
Let $S \times X \to X$ be a monoidal action of a monoid $S$ (with the identity $e$).
Assume that there exists a proper semitopological compactification
$\nu: X \hookrightarrow Y$ of $X$ which is $\{e\}$-topological
(that is, the action $S \times Y \to Y$ is continuous at every
$(e,y)$). If $F\subset X$ is a closed subset and $a\notin F$ then
there exist neighborhoods $U(e)$, $V(F)$ and $O(a)$ such that $UV
\cap UO =\emptyset$.
\end{lemma}
\begin{proof} Since $\nu: X \hookrightarrow Y$ is an embedding the closure
$cl(\nu(F))$ of $\nu(F)$ in $Y$ does not contain the point
$\nu(a)$. By the continuity of the action at every point $(e,y)$
(making use the Hausdorff axiom) it follows that for every $b \in
cl(\nu(F))$ there exist a neighborhood $U_b$ of $e$ and
neighborhoods $O_b$ of $\nu(a)$ and $V_b$ of $b$ in $Y$ such that
$U_bV_b \cap U_bO_b=\emptyset$. Now the standard compactness
argument easily completes the proof.

\end{proof}

Let $\pi: S \times X \to X$ be a jointly continuous semigroup
action. Up to an $S$-isomorphisms we can assume that $S$ and $X$
are disjoint sets. Denote by $S \sqcup_{\pi} X$ a new semigroup
defined as follows. As a set it is a {\it disjoint union} $S \cup
X$. The multiplication is defined by setting:

$a \circ b:=sx$ if $a=s \in S, \ b=x \in X$

$a \circ b:=s_1s_2$  if $a=s_1 \in S, \ b=s_2 \in S$

and

$a \circ b:=a$ in other cases.

Then $S \sqcup_{\pi} X$ is a topological semigroup which we call a
{\it $\pi$-generated semigroup}.


\begin{lemma} \label{l:constr}
Let $X$ be an $S$-space.
\ben
\item
The topological semigroup $P:=S \sqcup_{\pi} X$ is
 compactifiable (semi-compactifiable) if and only if $(S,X)$ is a
 compactifiable (resp.: semi-compactifiable) flow
and at the same time $S$ is a  compactifiable (resp.: semi-compactifiable) semigroup.
\item
The opposite topological semigroup $P^{op}:=(S \sqcup_{\pi}
X)^{op}$ is compactifiable if and only if $S^{opp}$ is a
compactifiable semigroup and the topology of $X$ admits a system
of $S$-contractive pseudometrics.

\een
\end{lemma}
\begin{proof}
(1): Observe that we have naturally defined equivariant inclusion of flows
$$(h,\a): (S,X) \rightrightarrows (P,P)=( S \sqcup_{\pi} X, S \sqcup_{\pi} X).$$
Therefore if $(P,P)$ is compactifiable then the same is true for
$(S,X)$ and $(S,S)$.

Conversely, every pair $\psi_1: S \hookrightarrow Y_1$ and
$\psi_2: X \hookrightarrow Y_2$ of proper $S$-compactifications
(one may assume that $Y_1$ and $Y_2$ are disjoint) defines a
proper $P$-compactification $\psi: P=S \sqcup_{\pi} X
\hookrightarrow Y_1 \sqcup Y_2$.

(2): If $P^{op}$ is compactifiable then $S^{op}$ being a
subsemigroup of $P^{op}$ is also compactifiable. Moreover, by
Theorem \ref{t:sem} there exists a system of right contractive
pseudometrics on $P^{op}=(S \sqcup_{\pi} X)^{op}$. Such a system
is clearly left contractive on $P$. It induces the desired system
of
$S$-contractive pseudometrics on $X$.

Conversely, suppose that $S^{op}$ is compactifiable and the
topology of $X$ is generated by a family $\Fcal_1:= \{d_i\}_{i \in
I}$ of $S$-contractive pseudometrics. By the first assumption and
Theorem \ref{t:sem} there exists a system $\Fcal_2:= \{\rho_j\}_{j
\in J}$ of left contractive pseudometrics on $S$. One can suppose
in addition that $d_i \leq 1$ and $\rho_j \leq 1$ for every $(i,j)
\in I \times J$.

Now define a new system $\Fcal_3=\Fcal_1 \cup \Fcal_2 \cup \{D\}$
on $P=S \sqcup_{\pi} X$ by setting $D(s,x)=D(x,s)=1$ for every $s
\in S, \ x \in X$ and $D(s_1,s_2)=D(x_1,x_2)=0$ for every $s_1,
s_2 \in S, \ x_1, \ x_2 \in X$. It is easy to verify that
$\Fcal_3$ is a system of left contractive pseudometrics on $P$
generating its topology. The same system is right contractive on
$P^{op}$. Hence by Theorem \ref{t:sem} we can conclude that
$P^{op}$ is compactifiable.
\end{proof}


\begin{examples} \label{ex:count}
Here we give some examples of noncompactifiable topological
semigroups and actions. \ben
\item
{\it The linear action of the compact multiplicative monoid
$S:=([0,1],\cdot)$ on $X:=[0,\infty)$ is not compactifiable}.
Moreover, every $f \in RUC_S(X)$ is necessarily constant.

 \sk Assuming the contrary let $f \in RUC_S(X)$ be nonconstant. Then $f(a)
 -f(b) =\eps >0$ for a pair $a, b \in X$. By definition of $RUC_S(X)$ there
 exists $\delta > 0$  such that $|f(u_1x) - f(u_2x)| < \eps$ for
 every triple $(u_1,u_2,x) \in U \times U \times X$, where
 $U:=[0,\delta)$. Choose $x_0 \in X$ such that $a <
 \delta x_0$ and $b < \delta x_0$. Take $u_1:=\frac{a}{x_0}$ and
 $u_2:=\frac{b}{x_0}$.
 Then $(u_1,u_2,x_0) \in U \times U \times X$ but $|f(u_1x_0) -f(u_2x_0)| =
 \eps$.

\sk

Note that in this example the acting monoid is a submonoid of
$\Theta(V)$ for $V:=\R$. As a corollary we get that the action
$\Theta(V) \times V \to V$ is not compactifiable for any nontrivial normed
space $V$.

\item
{\it The multiplicative monoid $S:=([0,\infty), \cdot)$ (and hence
also the multiplicative monoid $\R$ of all reals) is not
compactifiable. In fact the corresponding universal dynamical
compactification $S \to S^{RUC}$ is a singleton.}

\sk This follows directly from example (1).

Since $\Theta(V)^{op}$ is compactifiable and $\R$ is involutive
(even, commutative), as a corollary of our results we get that
$(\R,\cdot)$ is not embedded into $\Theta(V)$ for arbitrary normed
space $V$. As well as $(\R,\cdot)$ is not embedded as a
topological subsemigroup into $U:=C(I^{\omega}, I^{\omega})$.

\item  {\it The universal right topological semigroup
compactification $S \to S^{LMC}$ of $S:=([0,\infty), \cdot)$ is
injective but not proper (that is, $LMC(S)$ separates the points
but does not determine the original topology). Hence, $[0,\infty),
\cdot)$ is not semicompactifiable.}

\sk Let $M$ be the additive monoid $\R \cup \{\theta\}$ where
topologically $\theta$ is a point at $+\infty$ and algebraically
$\theta+x=x+\theta=\theta$ for every $x \in M$. In fact this
semigroup $M$ is a copy of the multiplicative semigroup
$[0,\infty)$ via the topological isomorphism $\R \cup \{\theta\}
\to [0,\infty), \ \a(\theta)=0, \a(x)=2^{-x}$ for all $x \in \R$.
Now note that by results of Hindman and Milnes \cite[chapter
5]{HM} the algebra $LMC(M)$ separates the points but does not
determine the original topology (see also the results of Section
\ref{s:fun}).

\item
{\it One-parameter additive semigroup action on a Polish phase
space which is not semi-compactifiable.} \sk This construction was
inspired by Ruppert \cite[Ch. II, Example 7.8]{Ru}. Let $\R_+
=([0,\infty), +)$ be the one parameter additive semigroup. Denote
by $[0, \infty]$ the Alexandrov compactification of $\R_+$. In the
product space $[0, \infty] \times [0, \infty]$ consider the
following subspace
$$X:= [0, \infty) \times [0, \infty) \cup \{(\infty, \infty)\}$$
Then $X$ is Polish being homeomorphic to a $G_{\delta}$-subset of
the 2-cell $[0,1] \times [0,1]$. Define now the desired continuous action by
$$
\pi: \R_+ \times X \to X, \hskip 0.3cm t(x,y)=(x,tx+y), \hskip 0.3cm
t(\infty, \infty)=(\infty, \infty)
$$

Define in $X$ the point $a:=(\infty, \infty)$ and the closed
subset $F:=[0,\infty) \times \{0\}$. Then for every neighborhood
$O(a)$ of $a$ and every neighborhood $U(0)$ of $0$ in $\R_+$ we
have $UF \cap O \neq \emptyset$. Now Lemma \ref{key} and Remark
\ref{r:comp}.2 imply that $X$ is not semi-compactifiable.

\item
{\it Compact monoid action on a discrete space which is not
semi-compactifiable.}
 \sk Let $S$ be the compact
monoid homeomorphic to the Cantor cube $C:=\{0,1\}^{\N_0}$ and
$X=\N_0:= \N \cup \{0\}$. Look at the semigroup $C$ as the product
semigroup of elementary 2-point multiplicative monoids $\{0,1\}$.
Define the desired action by
$$\pi: C \times \N_0 \to \N_0, \hskip 0.3cm  \pi(c,n)=c_nn,$$
where $c=(c_k)_{k \in \N_0} \in C$. In $\N_0$ choose the point
$a:=0$ and the closed subset $F:=\N$. Then for every neighborhood
$U({\bold 1})$ holds $a \in UF$ (where ${\bold 1}=(1,1,1,\cdots)$
is the identity of the monoid $C$). Hence Lemma \ref{key} and
Remark \ref{r:comp}.1 finish the proof.


%
%
%
%


\item {\it A topological semigroup $Q$ such that Q is compactifiable and the
opposite semigroup $Q^{op}$ is not semi-compactifiable.}
 \sk
We construct the desired semigroup as the $\pi$-generated
semigroup $P:=\{0,1\}^{\N_0} \sqcup_{\pi} \N_0$ for the flow
$(S,X)=(C,\N_0)$ described in (5). Then $P$ is not
semi-compactifiable by Lemma \ref{l:constr}.1. Then the opposite
semigroup $Q:=P^{op}$ is the desired one. Indeed, first of all
$Q^{op}=P$ is not semi-compactifiable.

Clearly, $S=\{0,1\}^{\N_0}$ is compactifiable being a compact
semigroup. Define the standard $0,1$ metric on the discrete space
$X:=\N_0$. Then this metric is contractive with respect to the
action of $S$ on $X$.
 By Lemma \ref{l:constr}.2 we conclude that $P^{op}=Q$
is compactifiable.


%



\item
{\it There exists a Banach space $V$ such that the monoid
$\Theta(V)$ 
is not semi-compactifiable.}
 \sk

Let $Q$ be the topological semigroup defined in (6). Then
$P:=Q^{op}$ is not semi-compactifiable. On the other hand $P$,
being the opposite semigroup of a compactifiable semigroup $Q$, is
a topological subsemigroup of $\Theta(V)$ for some Banach space
$V$ (see Theorem \ref{t:sem}). Therefore $\Theta(V)$ is not
semi-compactifiable, too.


%

\item
{\it Sorgenfrey line $(\R_s, +)$ is a non-compactifiable
topological monoid.}
 \sk
 This follows directly from Theorem \ref{para}.
 Moreover it is not hard to see that
$RUC(\R_s)=RUC(\R)$. That is, the universal dynamical
compactification $\R_s^{RUC}$ is just the greatest ambit
$\R^{RUC}$ (for the usual topological group $\R$ of the reals).

\item {\it For every Polish not locally compact topological group $G$
there exists a continuous action $\pi: G \times X \to X$ on a
Polish space $X$ such that the corresponding $\pi$-generated
Polish semigroup $P:=G \sqcup_{\pi} X$ is not
semi-compactifiable.}
 \sk By \cite{MeSc1} there exists a
 non-compactifiable Polish $G$-space $X$.
 Then the semigroup $P:=G \sqcup_{\pi} X$ is not semi-compactifiable.
 Indeed assuming the contrary it follows by Lemma \ref{l:constr}.1
 that $(G,X)$ is semi-compactifiable.
 Since $G$ is \v{C}ech-complete we get (see Remark \ref{r:comp}.1)
 that $X$ is $G$-compactifiable, a contradiction.


\een
\end{examples}

\sk

\sk \sk

\bibliographystyle{amsplain}

\end{document}